\documentclass[centertags,12pt]{amsart}
\usepackage{latexsym}
\usepackage{amssymb}

\newtheorem{theorem}{Theorem}[section]

\newtheorem{lemma}[theorem]{Lemma}
\newtheorem{corollary}[theorem]{Corollary}
\theoremstyle{definition}

\newtheorem{example}[theorem]{Example}
\newtheorem{remark}[theorem]{Remark}
\newtheorem*{remark*}{Acknowledgment}

\numberwithin{equation}{section}

\def\F{\mathbb F}
\def\Z{\mathbb Z}

\DeclareMathOperator{\N}{Norm}
\DeclareMathOperator{\ind}{ind}
\DeclareMathOperator{\tr}{tr}
\DeclareMathOperator{\Tr}{Tr}

\begin{document}
\title[Irreducible polynomials and Kloosterman sums]{Kloosterman sums, elliptic curves, and irreducible polynomials with prescribed trace and norm}

\author{Marko Moisio}
\email{mamo@uwasa.fi}
\address{Department of Mathematics and Statistics, University of Vaasa, P.O.
Box 700, FIN-65101 Vaasa, Finland}
\date{\today}

\begin{abstract}Let $\F_q$ ($q=p^r$) be a finite field. In this paper the number of irreducible polynomials of degree $m$ in $\F_q[x]$ with prescribed trace and norm coefficients is calculated in certain special cases and a general bound for that number is obtained improving the bound by Wan if $m$ is small compared to $q$. As a corollary, sharp bounds are obtained for the number of elements in $\F_{q^3}$ with prescribed trace and norm over $\F_q$ improving the estimates by Katz in this special case. Moreover, a characterization of Kloosterman sums over $\F_{2^r}$ divisible by three is given generalizing the earlier result by Charpin, Helleseth, and Zinoviev obtained only in the case $r$ odd. Finally, a new simple proof for the value distribution of a Kloosterman sum over the field $\F_{3^r}$, first proved by Katz and Livne, is given.  
\end{abstract}


\maketitle

\section{Introduction}

Let $\F_q$ be a finite field with $q = p^r$, and let $a,b\in\F_q$, $b\ne0$. Fairly little is known about the number $P_m(a,b)$ of irreducible polynomials $p(x)=x^m-ax^{m-1}+\cdots+(-1)^mb$ in $\F_q[x]$. Carlitz~\cite{Carlitz52} obtained the asymptotic formula 
\[
	P_m(a,b)=\frac{q^m-1}{mq(q-1)}+\mathcal O(q^{m/2})\quad (m\to\infty),
\]
and evaluated $\sum_b P_m(a,b)$ where $b$ runs over $\F_q^*$, and over the set of squares (resp. non-squares) in $\F_q^*$. Later Yucas~\cite{Yucas06} calculated elementarily the numbers $\sum_{a}P_m(a,b)$ and $\sum_{b}P_m(a,b)$ where $a,b$ run over $\F_q$. 

By the bijection $p(x)\mapsto (-1)^mb^{-1}x^mp(\frac 1x)$ we see that $P_m(a,b)$ equals the number of irreducible monic polynomials of degree $m$ in the arithmetic progression 
\[ 
	\{cx+d+f(x)x^2\mid f(x)\in\F_q[x]\} 
\]	
where $c=(-1)^{m+1}ab^{-1}$ and $d=(-1)^mb^{-1}$. Applying a general asymptotic bound on the number of primes on an arithmetic progression (see e.g.~\cite[p.40]{Rosen}) we actually have the asymptotic bound 
\[
P_m(a,b)=\frac{q^{m-1}}{m(q-1)}+\mathcal O\Big(\frac{q^{m/2}}m\Big)\quad (m\to\infty).  
\]
Finally, Wan~\cite[Thm. 5.1]{Wan} obtained the following effective bound: 
\begin{equation}\label{e:carwan}
	\left|P_m(a,b)-\frac{q^{m-1}}{m(q-1)}\right|\le\textstyle\frac 3mq^{\frac m2}. 
\end{equation}
For a more complete survey the reader is referred to~\cite{Cohen05}.
 
The bounds above are obtained by using Dirichlet L-series over $\F_q[x]$ and the Riemann's hypothesis for function fields over a finite field. In this paper we express $P_m(a,b)$ in terms of the numbers $N_t(a,b)$ of elements $x\in\F_{q^t}$ (with $t\mid m$) satisfying $Trace(x)=a$ and $Norm(x)=b$ ($Trace,Norm$ are from $\F_{q^m}$ onto $\F_q$), which, in turn, are expressed in terms of exponential sums. This opens up a possibility to calculate $P_m(a,b)$ explicitly in certain special cases. Moreover, we shall obtain an improvement of the bound~(\ref{e:carwan}) if $m$ is small compared to $q$, more precisely, if $m\le\frac32(q-1)$. If $a = 0$ the bound is obtained elementarily, but if $a\ne0$ this is done by linking the problem to the number of solutions of certain system of equations, and making use of the Katz bound~\cite{Katz}:
\begin{equation}\label{e:Katz}
\left|N_m(a,b)-\frac{q^m-1}{q(q-1)}\right|\le mq^{\frac{m-2}2}, 
\end{equation}
proved by using deep algebraic geometry.

The Katz bound with $m = 3$ plays a significant role in the proof by Cohen and Huczynska~\cite{Huc-Cohen} of the existence of a primitive free (normal) cubic polynomial with $a\,(\ne0)$ and $b$ fixed, which completed a general existence theorem (see also~\cite{Cohen00,Cohen-Huc,Cohen05}). We shall improve the Katz bound in this case. In fact, we get sharp lower and upper bounds for $N_3(a,b)$, and, as a corollary for $P_3(a,b)$, by using only the Hasse-Weil bound for elliptic curves together with a simple divisibility argument. 

Another special case where the Katz bound can be improved is the case $m = p^k$ for some $k$. Especially, if $p = 3$ (resp. $p = 2$) a result on the distribution of irreducible cubic (resp. quartic) polynomials in $\F_q[x]$ with trace and norm prescribed is obtained in terms of Kronecker class numbers by using the known value distribution of a Kloosterman sum over $\F_q$~\cite{Katz-Livne,La-Wo}.

Next, necessary and sufficient conditions for a Kloosterman sum over $\F_{2^r}$ divisible by three is given. In the case $r$ odd this result follows also from~\cite[Thm.~3]{chz}. Finally, a new proof for the value distribution of a Kloosterman sum over the field $\F_{3^r}$ is given. The proof uses only elementary properties of elliptic curves together with a result by Deuring~\cite{Deuring} which lies deeper: the knowledge of the number of isomorphism classes of elliptic curves over $\F_q$ having $q+1+t$ points with $\gcd(q,t)=1$.

\begin{remark*}The author is indebted to the anonymous referee for helpful comments and suggestions which improved the clarity of the paper considerably.  
\end{remark*}
 
\section{Basic formulae}

The aim of this section is to establish a link between the numbers $N_m(a,b)$ and $P_m(a,b)$, and to give basic formulae for $N_m(a,b)$ and $P_m(a,b)$ in terms of exponential sums. The formulae will be studied more closely in later sections.
\medskip

\begin{tabular}{cl}
$m,p,r$	&\ fixed positive integers, $m\ge2$, $p$ a prime\\
$\F_q$ &\ the finite field with $p^r$ elements\\
$a,b$ &\ fixed elements in $\F_q$, $b\ne0$\\
$P_m(a,b)$ &\ the number of irreducible polynomials\\
&\ $x^m-ax^{m-1}+\cdots+(-1)^mb\in\F_q[x]$\\
$t$ &\ a positive factor of $m$\\
$d$ &\ equals $\gcd(q-1,\frac mt)$\\
$\gamma_t$ &\ a primitive element of $\F_{q^t}$\\
$g$ &\ the primitive element of $\F_q$ defined by $g=\N_t(\gamma_t)$\\
$\tr_t(x)$ &\ the trace function from $\F_{q^t}$ onto $\F_q$\\
$\N_t(x)$ &\ the norm function from $\F_{q^t}$ onto $\F_q$\\
$S_t(a,b)$ &\ the set of the elements $x$ in $\F_{q^t}^*$ with\\
&\ $\tr_m(x)=a$ and $\N_m(x)=b$\\
$N_t(a,b)$ &\ the number of elements in $S_t(a,b)$\\
$\mu$ &\ the M\"obius function\\
$\chi$ and $e$ &\ the canonical additive characters of $\F_q$ and $\F_{q^t}$\\
$\mathcal X(\F_q)$ &\ the set of rational points on an algebraic\\
&\ curve $\mathcal X$ defined over $\F_q$
\end{tabular}
\bigskip

The following two lemmas relate the numbers $P_m(a,b)$ and $N_t(a,b)$:

\begin{lemma}\label{l:polyN}
\[
	P_m(a,b)=\frac1m\sum_{t\mid m}\mu(t)N_{m/t}(a,b).
\]
\end{lemma}

\begin{proof}Let
\[
	H_t(a,b)=\big|\{x\in\F_{q^t}^*\mid \tr_m(x)=a, \N_m(x)=b,\text{ and } x\not\in\F_{q^s}\text{ if } s<t\}\big|.
\]
Obviously $N_m(a,b)=\sum_{t\mid m}H_t(a,b)$, and now by M\"obius inversion formula 
\[
	H_m(a,b)=\sum_{t\mid m}\mu(t)N_{m/t}(a,b).
\]
But $H_m(a,b)=mP_m(a,b)$ completing the proof. 
\end{proof}
 
\begin{lemma}\label{l:bounds}Let $m=p_1^{e_1}\cdots p_k^{e_k}$ be the canonical prime number decomposition of $m$ ($p_1<p_2<\dots$), and let $m'=p_1\cdots p_k$. Then
\[
	N_m(a,b)-M_1m'/2\le mP_m(a,b) \le N_m(a,b)+M_2(m'-2)/2
\]
with $M_1=\max_{h}\{N_{m/h}(a,b)\}$, $M_2=\max_{s}\{N_{m/s}(a,b)\}$ where $h$ (resp. $s>1$) runs over the factors of $m'$ having odd (resp. even) number of prime factors. If $k=1$, set $M_2=0$.
\end{lemma} 

\begin{proof} Assume $k=1$. Now $mP_m(a,b)=N_{m}(a,b)-N_{m/p_1}(a,b)$ by Lemma~\ref{l:polyN}. Moreover, since $M_1=N_{m/p_1}(a,b)$ and $M_2=0$, the lemma follows in this case. Assume $k>1$. By Lemma~\ref{l:polyN} we have
\[
	mP_m(a,b)=N_m(a,b)+\sum_{s}N_{m/s}(a,b)-\sum_{h}N_{m/h}(a,b).
\]
Since $m'\ge 2^k$, we now get 
\begin{eqnarray*}
	mP_m(a,b)-N_m(a,b)&\ge&-M_1\sum_{h} 1 = -M_1\sum_{i=0}^{\lfloor k/2 \rfloor}\binom{k}{2i+1} = -M_12^{k-1}\\
	&\ge& -M_1m'/2.
\end{eqnarray*}
Moreover,
\begin{eqnarray*}
	mP_m(a,b)-N_m(a,b)&\le& M_2\sum_{s} 1 = M_2\sum_{i=1}^{\lfloor k/2 \rfloor}\binom{k}{2i} = M_2(2^{k-1}-1)\\
 &\le& M_2(m'-2)/2,
\end{eqnarray*}
and the proof is complete.
\end{proof}

Next we derive a formula for $N_t(a,b)$. First, we observe that if $a\ne0$, then  
\begin{equation}\label{e:x_in_S_t}
	x\in S_t(a,b)\ \Leftrightarrow\ p\nmid\textstyle\frac mt\text{ and } \tr_t(x)=\frac tm a 
	\text{ and } \N_t(x^\frac mt)=b, 
\end{equation}
and if $a=0$, then
\begin{align}\label{e:x_in_S_t_a=0}
x\in S_t(a,b)\ \Leftrightarrow\ &p\mid\textstyle\frac mt \text{ and } \N_t(x^\frac mt)=b,\\ 
&\text{ or, }p\nmid\frac mt,\ \tr_t(x)=0, \text{ and } \N_t(x^\frac mt)=b.\notag
\end{align} 

Second, we see that
\begin{align}\label{e:congruence}
	\N_t(x^\frac mt)=b &\Leftrightarrow \textstyle\frac mt i\equiv\ind_gb\pmod{q-1}\\
	&\Leftrightarrow d\mid\ind_gb\text{ and }i=i_0+\textstyle\frac{q-1}dj\notag
\end{align}
where $j$ runs over the set $\{0,\dots,\frac{q^t-1}{(q-1)/d}-1\}$, and $i_0$ is a solution of the congruence $\frac{m}{dt}i\equiv\frac{\ind_gb}d\pmod{\frac{q-1}d}$. 

\begin{lemma}\label{l:N_t} Assume $p\nmid\frac mt$ and $d\mid\ind_g b$. Let $i_0$ be a solution of the congruence $\frac{m}{dt}i\equiv\frac{\ind_gb}d\pmod{\frac{q-1}d}$ and let $a_0=\frac tm a$. 
Then
\[
	N_t(a,b)=\frac d{q(q-1)}(q^t-1+\sigma_t(a,b)),
\] 
where
\begin{equation}\label{e:def_S}
	\sigma_t(a,b)=\sum_{c\in\F_q^*}\chi(-ca_0)\sum_{x\in\F_{q^t}^*}e(c\gamma_t^{i_0}x^{\frac{q-1}d}).
\end{equation}
\end{lemma}

\begin{proof} Let $\alpha$ be an element in $\F_{q^t}$ with $\tr_t(\alpha)=t/m$. Now, by~\eqref{e:congruence} and by the orthogonality of characters we get
\begin{align*}
qN_t(a,b)&=\sum_{j=0}^{\frac{q^t-1}{(q-1)/d}-1}\sum_{c\in\F_q}
\chi(c\tr_m(\gamma_t^{i_0+\frac{q-1}dj}-\alpha a))\\
&=\sum_{c\in\F_q}\chi(-ca)\sum_{j=0}^{\frac{q^t-1}{(q-1)/d}-1}e({\textstyle\frac mt} c\gamma_t^{i_0}
\gamma_t^{\frac{q-1}dj})\\
&\stackrel{c\mapsto \frac tm c}{=}\frac d{q-1}\sum_{c\in\F_q}\chi(-ca_0) \sum_{x\in\F_{q^t}^*}e(c\gamma_t^{i_0}x^{\frac{q-1}d})\\
&=\frac d{q-1}(q^t-1+\sigma_t(a,b)).
\end{align*}
\end{proof}

\section{Zero trace}

In this section we assume that $a=0$ and simplify formula~\eqref{e:def_S} by using Gauss sums and some very elementary group theory. This will enable us to obtain an improvement of the Katz bound and the Wan bound in the case $a=0$. We use the following notations:
\medskip

\begin{tabular}{cl}
$H_n$ &\ the subgroup of order $n$ of the multiplicative\\
&\  character group of $\F_q$\\
$\lambda_0$ &\ the trivial character of $H_n$\\
\end{tabular}
\bigskip

For a multiplicative character $\psi$ of $\F_{q^t}$ we define a Gauss sum
\[
	G(\psi):=\sum_{x\in\F_{q^t}^*}e(x)\psi(x).
\]
 
\begin{lemma}\label{l:gaussreps}Let $n$ be a factor of $q-1$ and let $\alpha\in\F_{q^t}^*$. Then
\[
	\sum_{x\in\F_{q^t}^*}e(\alpha x^n)=
	\sum_{\lambda\in H_n}G(\bar\lambda\circ\N_t)\lambda(\N_t(\alpha)),
\]
where $\bar\lambda=\lambda^{-1}$.
\end{lemma}

\begin{proof}It is easy to see~\cite[p.217]{Lidl97} that 
\[
	\sum_{x\in\F_{q^t}^*}e(\alpha x^n)=
	\sum_{\psi\in H_n^\prime}G(\bar\psi)\psi(\alpha),
\]
where $H_n^\prime$ is the subgroup of order $n$ of the multiplicative character group of $\F_{q^t}$. But the surjectivity of $\N_t$ implies that $H_n^\prime=\{\lambda\circ\N_t\mid \lambda\in H_n\}$. 
\end{proof}

Assume $a=0$, $p\nmid\frac mt$, and $d\mid\ind_g b$. Now, by Lemma~\ref{l:gaussreps} we get
\begin{align*}
\sigma_t(0,b)&=\sum_{c\in\F_q^*}\sum_{x\in\F_{q^t}^*}e(c\gamma_t^{i_0}x^{\frac{q-1}d})\\
&=\sum_{c\in\F_q^*}\sum_{\lambda\in H_{\frac{q-1}d}}G(\bar\lambda\circ\N_t)
\lambda(c^tg^{i_0})\\
&=\sum_{\lambda\in H_{\frac{q-1}d}}G(\bar\lambda\circ\N_t)\lambda(g^{i_0})
\sum_{c\in\F_q^*}\lambda(c^n),
\end{align*}
where $n=\gcd(q-1,t)$. Since
\[
	\sum_{c\in\F_q^*}\lambda(c^n)=\sum_{c\in\F_q^*}\lambda^n(c)=
	\begin{cases}
	0&\text{ if }\lambda^n\ne\lambda_0,\\
	q-1&\text{ if }\lambda^n=\lambda_0\Leftrightarrow\text{ if }\lambda\in H_n\cap H_{\frac{q-1}d}, 
	\end{cases}
\]
we get
\[
	\sigma_t(0,b)=(q-1)\sum_{\lambda\in H_s}G(\bar\lambda\circ\N_t)\lambda(g^{i_0})
	=(q-1)\sum_{x\in\F_{q^t}^*}e(\gamma_t^{i_0}x^s),
\]
where $s=\gcd(n,\frac{q-1}d)$. Thus,
\[
	N_t(0,b)=\frac dq\Big(\frac{q^t-1}{q-1}+\sum_{x\in\F_{q^t}^*}e(\gamma_t^{i_0}x^{s})\Big),
\]
implying the following

\begin{theorem}\label{t:N_bound}Assume $p\nmid\frac mt$ and $d\mid\ind_g b$. Then,
\[
	N_t(0,b)=d\Big(\frac{q^{t-1}-1}{q-1}+\frac1q\sum_{x\in\F_{q^t}}e(\gamma_t^{i_0}x^{s})\Big),
\]
where $s=\gcd(t,\frac{q-1}d)$ and $d=\gcd(\frac mt,q-1)$. 
\end{theorem}

Theorem~\ref{t:N_bound} and the Weil bound (see e.g.~\cite[p.223]{Lidl97}) imply an improvement of the Katz-bound (see~\eqref{e:Katz}) in the case $a=0$:

\begin{corollary}\label{c:bound_a_0}
\[
	\left|N_m(0,b)-\frac{q^{m-1}-1}{q-1}\right|\le (s-1)q^{\frac{m-2}2},
\] 
where $s=\gcd(m,q-1)$.
\end{corollary}

We can now prove an improvement of the Wan bound (see~\eqref{e:carwan}) in the case $a=0$ and $m\le\frac32(q-1)$:
\begin{corollary}
\[
	\left|P_m(0,b)-\frac{q^{m-1}-1}{m(q-1)}\right|\le\frac{s-1}m q^{\frac{m-2}2}+\frac{q^{\frac m2}-1}{q-1}
		<\frac2{q-1}q^{\frac m2}, 
\]
where $s=\gcd(m,q-1)$.
\end{corollary}

\begin{proof}Since $d\le m/t$, it follows from~\eqref{e:x_in_S_t_a=0} and~\eqref{e:congruence} that the numbers $M_1$ and $M_2$ in Lemma~\ref{l:bounds} satisfy 
\[
	M_2<M_1\le p_1\frac{q^{m/p_1}-1}{q-1}\le2\frac{q^{m/2}-1}{q-1},
\]
and now, by Lemma~\ref{l:bounds} and Corollary~\ref{c:bound_a_0}, we get 
\[
	\left|mP_m(0,b)-\frac{q^{m-1}-1}{q-1}\right|\le(s-1)q^{\frac{m-2}2}+m\frac{q^{m/2}-1}{q-1}.
\]
\end{proof}

By Lemma~\ref{l:polyN}, Theorem~\ref{t:N_bound},~\eqref{e:congruence}, and~\eqref{e:x_in_S_t_a=0} we get explicit expressions for $P_m(0,b)$ e.g. in the following special cases:

\begin{example}If $\gcd(p,m,q-1)=1$, then
\[
	P_m(0,b)=\frac1{m(q-1)}\sum_{t\mid m}\mu\Big(\frac mt\Big)(q^{t-1}-1).
\]
\end{example}

\begin{example}\label{ex:m_p_k} If $m=p^k>1$, then
\[
	mP_m(0,b)=\frac{q^{m-1}-1}{q-1}-\frac{q^{m/p}-1}{q-1}.
\]
\end{example}

\section{Non-zero trace}

In this section we assume that $a\ne0$. This case is much harder than the zero trace case, and we are not able to find such a simple expression for $N_t(a,b)$ as in case $a=0$. The best we can do is to give $N_t(a,b)$ in terms of the number of solutions of a system of equations, and estimate that number by using the Katz-bound. This method will lead us to get an improvement of the Wan bound also in the case $a\ne0$.  
 
Let $n=(q-1)/d$. By Lemma~\ref{l:gaussreps} and by substitution $c\mapsto -a_0^{-1}c$ we see that 
$\sigma_t(a,b)$ (see~\eqref{e:def_S}) can be written in the form
\begin{align*}
\sigma_t(a,b)&=\sum_{c\in\F_q^*}\chi(c)\sum_{\lambda\in H_n}
G(\bar\lambda\circ\N)\lambda\big(c^t(-a_0)^{-t}g^{i_0}\big)\\
&=\sum_{\lambda\in H_n}G(\bar\lambda\circ\N)\lambda\big(g^{i_0}(-a_0)^{-t}\big)
\sum_{c\in\F_q^*}\chi(c)\lambda^t(c)\\
&=\sum_{\lambda\in H_n}G(\bar\lambda\circ\N)G(\lambda^t)
\lambda\big(g^{i_0}(-a_0)^{-t}\big).
\end{align*}
Let $c=g^{i_0}a_0^{-t}$ and use Davenport-Hasse Theorem~\cite[p.197]{Lidl97} to get
\begin{align*}
\sigma_t(a,b)&=(-1)^{t-1}\sum_{\lambda\in H_n}G(\bar\lambda)^tG(\lambda^t)\lambda\big((-1)^tc\big)\\
&=(-1)^{t-1}\sum_{\lambda\in H_n}G(\lambda)^tG(\bar\lambda^t)\bar\lambda\big((-1)^tc\big).
\end{align*}

Now, by the definition of a Gauss sum we get 
\[
	\sigma_t(a,b)=(-1)^{t-1}\sum_{x_1,\dots,x_t,u\in\F_q^*}\chi(x_1+\cdots+x_t+u)
	\sum_{\lambda\in H_{n}}\lambda\big(x_1\cdots x_t(-u)^{-t}c^{-1}\big),
\]
and consequently, by substituting $x_1\mapsto -ux_1,\dots,x_t\mapsto -ux_t$, we obtain
\begin{align}\label{e:sigma_a_0}
\sigma_t(a,b)=(-1)^{t-1}&\sum_{x_1,\dots,x_t,u\in\F_q^*}\chi(-u(x_1+\cdots+x_t-1))\\
&{\hskip15pt}\times\sum_{\lambda\in H_n}\lambda(x_1\cdots x_tc^{-1}).\notag
\end{align}

We can now prove the following

\begin{theorem}\label{t:N_t_pair}If $a\ne0$, $p\nmid\frac mt$, and $d\mid\ind_g b$, then
\[
	N_t(a,b)=\frac{d(q^t-1)}{q(q-1)}+(-1)^{t-1}\Big(\sum_{i=0}^{d-1}N(c_i)-\frac{d(q-1)^t}{q(q-1)}\Big),
\]
where $N(c_i)$ is the number of solutions of
\[
\begin{cases}
	x_1+\cdots+x_t&=1\\
	x_1\cdots x_t&=c_i\\
\end{cases}
\]
in $\F_q^t$ with $c_i=g^{\frac{q-1}d i + i_0}a_0^{-t}$.
\end{theorem}

\begin{proof}Let $n=(q-1)/d$ and $c=g^{i_0}a_0^{-t}$. The orthogonality of characters implies that
\[
	q(q-1)N(c_i)=\sum_{x_1,\dots, x_t\in\F_q^*}\sum_{u\in\F_q}\chi(u(x_1+\cdots+x_t-1))
	\sum_{\lambda\in H_{q-1}}\lambda(c_i^{-1}x_1\cdots x_t),
\] 
and consequently
\begin{align*}
	q(q-1)\sum_{i=0}^{d-1}N(c_i)=&\sum_{x_1,\dots, x_t\in\F_q^*}\sum_{u\in\F_q}\chi(u(x_1+\cdots+x_t-1))\\
	&\times\sum_{\lambda\in H_{q-1}}\sum_{i=0}^{d-1}\lambda(c_i^{-1}x_1\cdots x_t).
\end{align*}
Here
\begin{align*}
	\sum_{i=0}^{d-1}\lambda(c_i^{-1}x_1\cdots x_t)&=\lambda(c^{-1}x_1\cdots x_t)\sum_{i=0}^{d-1}\lambda(g^{-ni})\\
	&=\begin{cases}
	\lambda(c^{-1}x_1\cdots x_t)d&\text{if }\lambda\in H_n,\\
	0,&\text{otherwise},
	\end{cases}
\end{align*}
and now, by~\eqref{e:sigma_a_0}, we get
\begin{align*}
	q(q-1)\sum_{i=0}^{d-1}N(c_i)
	&=d\sum_{x_1,\dots, x_t\in\F_q^*}\sum_{u\in\F_q}\chi(u(x_1+\cdots+x_t-1))\\
	&{\hskip30pt}\times\sum_{\lambda\in H_n}\lambda(c^{-1}x_1\cdots x_t)\\
	&=d(-1)^{t-1}\sigma_t(a,b)+d\sum_{x_1,\dots, x_t\in\F_q^*}\sum_{\lambda\in H_n}\lambda(c^{-1}x_1\cdots x_t).
\end{align*}
Here
\[
	\sum_{x_1,\dots, x_t\in\F_q^*}\sum_{\lambda\in H_n}\lambda(c^{-1}x_1\cdots x_t)
	=\sum_{\lambda\in H_n}\Big(\sum_{x\in\F_q^*}\lambda(c^{-1}x)\Big)^t=(q-1)^t
\]
and it follows that
\[
	\sigma_t(a,b)=(-1)^{t-1}\Big(\frac{q(q-1)}d\sum_{i=0}^{d-1}N(c_i)-(q-1)^t\Big).
\]

Lemma~\ref{l:N_t} now completes the proof.
\end{proof}

\begin{lemma}\label{l:key}Let $n$ be a positive integer and let $c\in\F_q^*$. The number $N(c)$ of solutions $(x_1,\dots,x_n)$ in $\F_q^n$ of
\[
\begin{cases}
	x_1+\cdots+x_n&=1\\
	x_1\cdots x_n&=c\\
\end{cases}
\]
satisfies
\[
	\Bigg|N(c)-\frac{(q-1)^n}{q(q-1)}\Bigg|\le nq^{\frac{n-2}2}.
\]
\end{lemma}

\begin{proof}Choose $m=t=n$, and $a=1$, $b=c$. Now $d=\gcd(\frac mt,q-1)=1$ and we choose $i_0=\ind_g b$ (see~(\ref{e:congruence})). Now $c=g^{i_0}/a^t$, and by Theorem~\ref{t:N_t_pair} we get
\[
	N_n(a,b)=\frac{q^n-1}{q(q-1)}+N(c)-\frac{(q-1)^n}{q(q-1)}
\]
or, equivalently,
\[
	N(c)-\frac{(q-1)^n}{q(q-1)}=N_n(a,b)-\frac{q^n-1}{q(q-1)}.
\]
The Katz bound~\eqref{e:Katz} now completes the proof.
\end{proof}

We are now able to prove an improvement of the Wan bound~\eqref{e:carwan} in the case $a\ne0$ and $m\le\frac32(q-1)$:

\begin{corollary}\label{c:a_ne_0}Let $a,b\in\F_q^*$. Then,
\[
	\left|P_m(a,b)-\frac{q^m-1}{mq(q-1)}\right|
	\le q^{\frac{m-2}2}+\frac{q^{\frac m2}-1}{q(q-1)}+{\textstyle\frac m2}q^{\frac{m-4}4}
	<\frac2{q-1}q^{\frac m2}.
\]
\end{corollary}
 
\begin{proof} If $p\mid\frac mt$ or $d\nmid\ind_gb$, then $N_t(a,b)=0$ by~\eqref{e:x_in_S_t} and~\eqref{e:congruence}. Assume $p\nmid\frac mt$ and $d\mid\ind_gb$. If $t$ is even Theorem~\ref{t:N_t_pair} implies
\[
	N_t(a,b)\le\frac{d(q^t-1)}{q(q-1)}-dN(c)+\frac{d(q-1)^t}{q(q-1)}
\]
for some $c\in\F_q^*$. Now, by Lemma~\ref{l:key}
\[
	N_t(a,b)\le\frac{d(q^t-1)}{q(q-1)}-d\Big(\frac{(q-1)^t}{q(q-1)}-tq^{\frac{t-2}2}\Big)+\frac{d(q-1)^t}{q(q-1)}
	=\frac{d(q^t-1)}{q(q-1)}+dtq^{\frac{t-2}2}.
\]
Since $d\le m/t$, we get
\begin{equation}\label{e:N_t_bound}
	N_t(a,b)\le\frac{m(q^t-1)}{tq(q-1)}+mq^{\frac{t-2}2}.
\end{equation}

If $t$ is odd, then
\[
	N_t(a,b)\le\frac{d(q^t-1)}{q(q-1)}+dN(c)-\frac{d(q-1)^t}{q(q-1)}
\]
for some $c\in\F_q^*$, and
\[
	N_t(a,b)\le\frac{d(q^t-1)}{q(q-1)}+d\Big(\frac{(q-1)^t}{q(q-1)}+tq^{\frac{t-2}2}\Big)-\frac{d(q-1)^t}{q(q-1)}
	=\frac{d(q^t-1)}{q(q-1)}+dtq^{\frac{t-2}2}.
\]
Hence, bound~\eqref{e:N_t_bound} holds in this case too. 

Now, by~\eqref{e:N_t_bound}, it is clear that the numbers $M_1$ and $M_2$ in Lemma~\ref{l:bounds} satisfy
\[
	M_2<M_1\le2\frac{q^{\frac m2}-1}{q(q-1)}+mq^{\frac{m-4}4},
\]
and consequently, by Lemma~\ref{l:bounds} and by the Katz-bound~\eqref{e:Katz}, we get 
\[
  \left|mP_m(a,b)-\frac{q^m-1}{q(q-1)}\right|<mq^{\frac{m-2}2}+m\frac{q^{\frac m2}-1}{q(q-1)}+\textstyle\frac{m^2}2q^{\frac{m-4}4}.
\]
Hence,
\begin{align*}
\left|P_m(a,b)-\frac{q^m-1}{mq(q-1)}\right|&\le q^{\frac{m-2}2}+\frac{q^{\frac m2}-1}{q(q-1)}
+{\textstyle\frac m2}q^{\frac{m-4}4}\\
&<\frac1q\Big(1+\frac1{q-1}+{\textstyle\frac m2}q^{-\frac m4q}\Big)q^{\frac m2}\\
&=\Big(\frac1{q-1}+{\textstyle\frac m2}q^{-\frac{m+4}4}\Big)q^{\frac m2}.
\end{align*}	
Obviously $\frac m2q^{-\frac{m+4}4}<1/(q-1)$, and so the proof is complete.
\end{proof}
 
\section{Cubics and cubic extensions}

In this section we assume that $m=3$. Now the system of equations defined in the previous section is of degree 3, and therefore we can give $N_3(a,b)$, and also $P_3(a,b)$, in terms of the number of rational points on a cubic curve defined over $\F_q$. Some elementary manipulations of cubic curves together with the Hasse-Weil bound for elliptic curves, and the link between $P_3(a,b)$ and $N_3(a,b)$, will then lead to a sharp bound for $N_3(a,b)$, which is also an improvement of the Katz bound in the case $m=3$. The following result is a key for such an improvement:

\begin{theorem}\label{t:P_3_a_ne0}Let $c=ba^{-3}$, and let $\mathcal X$ be the projective curve over $\F_q$ 
defined by
\[
	\mathcal X : y^2+cy+xy=x^3.
\]
Then, $N_3(a,b)=|\mathcal X(\F_q)|$ and
\[
	P_3(a,b)=\textstyle\frac13(|\mathcal X(\F_q)|-\epsilon),
\]
where
\[
	\epsilon=
	\begin{cases}
	1&\text{if $p\ne3$ and $c=\frac 1{27}$},\\
	0&\text{otherwise}.
	\end{cases}
\]
\end{theorem}

\begin{proof}Let $m=3$ and apply Theorem~\ref{t:N_t_pair} with $t=3$ to get
\begin{equation}\label{e:N_3}
	N_3(a,b)=\frac{q^3-1}{q(q-1)}+N(c_0)-\frac{(q-1)^3}{q(q-1)}=N(c_0)+3,
\end{equation}
where $N(c_0)$ is the number of solutions $(x,y,z)$ in $\F_q^3$ of 
\[
\begin{cases}
	x+y+z&=1\\
	xyz&=c_0\\
\end{cases}
\]
with $c_0=g^{i_0}/a^3=b/a^3$, or, equivalently, $N(c_0)$ is the number of solutions of
\begin{equation}\label{e:ell_eq}
	x^2y+xy^2-xy=-c_0,
\end{equation}
in $\F_q^3$. 

Equation~(\ref{e:ell_eq}) defines an affine component of the projective curve defined by
\[
	\mathcal X':y^2+y-xy=-c_0x^3,
\] 
and that affine component has exactly three points at infinity. Hence 
\[
	N_3(a,b)=|\mathcal X'(\F_q)|
\] 
by~\eqref{e:N_3}. By multiplying both sides of the equation of $\mathcal X'$ by $c_0^2$ and then substituting 
$x\mapsto -c_0^{-1}x$ and $y\mapsto c_0^{-1}y$ we see that $\mathcal X'$ is isomorphic over $\F_q$ to $\mathcal X$.

It follows from Lemma~\ref{l:polyN} that $3P_3(a,b)=N_3(a,b)-N_1(a,b)$, and by~\eqref{e:x_in_S_t} $x\in S_1(a,b)$ if and only if $p\ne3$ and $b=(a/3)^3$. This completes the proof.  
\end{proof}

\begin{corollary}\label{c:singular}Assume $p\ne3$, and let $b=(a/3)^3$. Then
\[
	P_3(a,b)=\textstyle\frac13(q\pm 1),
\]
where the sign is plus if $p\equiv2\pmod3$ and $2\nmid r$, and otherwise the sign is minus.
\end{corollary}

\begin{proof}If $p=2$ then the equation of $\mathcal X$ is $y^2+(x+1)y=x^3$. This equation has only three solutions 
with $x=0,1$. By substituting $y\mapsto (x+1)y$, we can write the equation in the form $y^2+y=x^3/(x+1)^2$, and then by substitutintg $x\mapsto x+1$ we get the equation
\begin{equation}\label{e:m=3_p=2}
 y^2+y=x+1+x^{-1}+x^{-2}.
\end{equation} 
Since the absolute trace of $x^{-1}+x^{-2}$ equals zero we have $\chi(x^{-1}+x^{-2})=1$, and therefore equation~\eqref{e:m=3_p=2} has exactly
\begin{align*}
	\sum_{x\in\F_q^*\setminus\{1\}}(1+\chi(x+1))=q-2+\chi(1)\Big(\sum_{x\in\F_q^*}\chi(x)-\chi(1)\Big)
	=q-3-\chi(1)
\end{align*}
solutions in $\F_q^2$ with $x\ne0,1$. Hence, in the case $p=2$, $|\mathcal X(\F_q)|=q-3-\chi(1)+3+1$ and $P_3(a,b)=(q-\chi(1))/3$.

Assume $p\ne2$ and write the equation $y^2+cy+xy=x^3$ in the form $(y+\frac12(c+x))^2=x^3+\frac14(c+x)^2$. 
Substitute $y\mapsto y-\frac12(c+x)$ to get 
\[
	y^2=x^3+\textstyle\frac14x^2+\frac12cx+\frac{c^2}4=(x+\frac19)^2(x+\frac1{36}).  
\]
Finally, by substituting $x\mapsto x-\frac19$, we see that $\mathcal X$ is isomorphic over $\F_q$ to
\[
	C : y^2=x^2(x-\textstyle\frac1{12}).
\]

Let $F$ be the set of finite points of $C$ and let $F'$ be the set of finite points of the curve $C'$ defined over $\F_q$ by
\[
	C' : z^2=u-\textstyle\frac1{12}.
\]
We note that the map $(x,y)\mapsto (u=x,z=\frac yx)$ from $F\setminus\{(0,0)\}$ to $F'$ is injective, and it follows that $|F|=|F|'\pm 1$ depending on whether the equation $z^2=-\frac1{12}$ has, or has not, a solution in $\F_q$. Hence, $|C(\F_q)|=|C'(\F_q)|\pm 1=q+1\pm1$, and now, by Theorem~\ref{t:P_3_a_ne0}, we get $P_3(a,b)=\frac13(q+1\pm1-1)$.
\end{proof}

We can now prove an improvement of the Katz bound in the case $m=3$: 

\begin{theorem}\label{t:N_3_floor_ceil}Let $a,b\in\F_q$, $b\ne0$. Then 
\[
	3\left\lceil\frac{q+1-2\sqrt q}3\right\rceil\le N_3(a,b)\le3\left\lfloor\frac{q+1+2\sqrt q}3\right\rfloor.
\]
\end{theorem}

\begin{proof}By Lemma~\ref{l:polyN} we have 
\[
	3P_3(a,b)=N_3(a,b)-N_1(a,b).
\]
Assume first that $a=0$. If $p=3$ then $N_1(a,b)=1$ by~\eqref{e:x_in_S_t_a=0}, and $3P_3(a,b)=q+1-1$ by Example~\ref{ex:m_p_k}. Hence, $N_3(a,b)=q+1$, and the theorem follows in the case $a=0$ and $p=3$. 

If $p\ne 3$ then $N_1(a,b)=0$ by~\eqref{e:x_in_S_t_a=0}, and Corollary~\ref{c:bound_a_0} now implies
\begin{equation}\label{e:P_3_bound}
	q+1-2\sqrt q\le 3P_3(a,b)\le q+1+2\sqrt q, 
\end{equation}
and therefore 
\[
	\left\lceil\frac{q+1-2\sqrt q}3\right\rceil\le P_3(a,b)\le\left\lfloor\frac{q+1+2\sqrt q}3\right\rfloor.
\]
Since $3P(a,b)=N_3(a,b)$, the proof is complete in case $a=0$ and $p\ne3$.

Assume next that $a\ne0$. It is easy to see that if $\mathcal X:y^2+cy+xy=x^3$ is singular then $p\ne3$ and $c=\frac1{27}$. Hence, if $\mathcal X$ is singular then $q\pm1=3P_3(a,b)=N_3(a,b)-1$ by Corollary~\ref{c:singular} and Theorem~\ref{t:P_3_a_ne0}, and therefore $N_3(a,b)=q\pm1+1$ proving the theorem if $\mathcal X$ is singular. 

Assume that $\mathcal X$ is non-singular.  Now, by the proof of Corollary~\ref{c:singular}, we see that $p=3$ or $c\ne\frac1{27}$, and therefore $3P_3(a,b)=|\mathcal X(\F_q)|=N_3(a,b)$ by Theorem~\ref{t:P_3_a_ne0}. Now, since $\mathcal X$ elliptic, the Hasse-Weil bound (see e.g.~\cite[p.91]{Washington}) implies that the bounds in~\eqref{e:P_3_bound} hold in this case too, and the proof is complete.
\end{proof}

\begin{remark}The bounds in Theorem~\ref{t:N_3_floor_ceil} are sharp. Take $q=5$, for example. If $a=b=1$, we have $N_3(a,b)=|\mathcal X(\F_q)|=9=3\lfloor(5+1+2\sqrt 5)/3\rfloor$. If $a=1$, $b=2$, we 
have $N_3(a,b)=|\mathcal X(\F_q)|=3=3\lceil(5+1-2\sqrt 5)/3\rceil$. These calculations can be verified e.g. by MAGMA.
\end{remark}

\section{Degree a power of the characteristic}\label{power}

An improvement of the Katz bound can also be obtained in the special case $m=p^k\,(>2)$ as we shall see in this section. The key point is that in this case the number of solutions of our system of equations, and therefore $N_m(a,b)$ and $P_m(a,b)$, can be given in terms of hyper-Kloosterman sums over $\F_q$ which can be estimated by the Deligne bound obtained in~\cite{Deligne} (see also~\cite[p.254]{Lidl97}). 

In the special cases $(p,m)=(3,3),(2,4)$ we can go even further since then we can use the known value distributions of Kloosterman sums to get fairly precise information on the distribution of the irreducible cubic and quartic polynomials over the fields $\F_{3^r}$ and $\F_{2^r}$, respectively. These cases are condidered in subsections~\ref{cubics} and~\ref{quartics}.

For a positive integer $n$ and $c$ in $\F_q^*$ let $k_n(c)$ be an $n$-dimensional Kloosterman sum (or a hyper-Kloosterman sum): 
\[
	k_n(c)=\sum_{x_1,\dots,x_n\in\F_q^*}\chi\left(x_1+\cdots+x_n+\frac{c}{x_1\cdots x_n}\right).
\]

\begin{theorem}\label{t:N_p_kl}Assume $m=p^k>2$, and let $a,b\in\F_q^*$. Then,
\[
	N_m(a,b)=\frac{q^{m-1}-1}{q-1}+(-1)^{m-1}k_{m-2}(c),
\]
where $c=b/a^m$. Moreover,
\[
	\left|N_m(a,b)-\frac{q^{m-1}-1}{q-1}\right|\le(m-1)q^{\frac{m-2}2}.
\]
\end{theorem}

\begin{proof}Apply Theorem~\ref{t:N_t_pair} with $m=t$ to get
\begin{equation}\label{e:N_p}
N_m(a,b)=\frac{q^m-1}{q(q-1)}+(-1)^{m-1}\left(N(c)-\frac{(q-1)^m}{q(q-1)}\right),
\end{equation}
where $N(c)$ is the number of solutions of
\[
\begin{cases}
	x_1+\cdots+x_m&=1\\
	x_1\cdots x_m&=c
\end{cases}
\]

Obviously $N(c)$ is equal to the number of solutions of
\[
x_1+\cdots+x_{m-1}+\frac{c}{x_1\cdots x_{m-1}}-1=0,
\]
and therefore, by the orthogonality of characters, we get 
\begin{align*}
qN(c)&=\sum_{x_1,\dots,x_{m-1}\in\F_q^*}\sum_{u\in\F_q}\chi\left(u\Big(x_1+\cdots+x_{m-1}
+\frac c{x_1\cdots x_{m-1}}-1\Big)\right)\\
&=\sum_{u\in\F_q^*}\chi(-u)\sum_{x_1,\dots,x_{m-1}\in\F_q^*}
\chi\left(ux_1+\cdots+ux_{m-1}+\frac{uc}{x_1\cdots x_{m-1}}\right)\\
&{\hskip10pt}+(q-1)^{m-1}.
\end{align*}
Now, by substitutions $x_1\mapsto x_1/u,\dots,x_{m-1}\mapsto x_{m-1}/u$, and by noting that $x\mapsto x^m$ is a permutation of $\F_q$, we get
\begin{align*}
&qN(c)-(q-1)^{m-1}\\
&=\sum_{u\in\F_q^*}\chi(-u)\sum_{x_1,\dots,x_{m-1}\in\F_q^*}\chi\left(x_1+\cdots+x_{m-1}+
\frac{u^mc}{x_1\cdots x_{m-1}}\right)\\
&=\sum_{u\in\F_q^*}\chi(-u)\sum_{x_1,\dots,x_{m-1}\in\F_q^*}
\chi\left(x_1^m+\cdots+x_{m-1}^m+\frac{u^mc}{x_1^m\cdots x_{m-1}^m}\right)\\
&=\sum_{u\in\F_q^*}
\chi(-u)\sum_{x_1,\dots,x_{m-1}\in\F_q^*}\chi\left(\Big(x_1+\cdots+x_{m-1}+\frac{uc^{1/m}}
{x_1\cdots x_{m-1}}\Big)^m\right)\\
&=\sum_{x_1,\dots,x_{m-1}\in\F_q^*}\chi(x_1+\cdots+x_{m-1}) \sum_{u\in\F_q^*}\chi\left(u\Big(\frac{c^{1/m}}{x_1\cdots x_{m-1}}-1\Big)\right).
\end{align*}
The inner sum equals $q-1$ or $-1$ according as $x_1\cdots x_{m-1}$ is, or is not, equal to $c^{1/m}$. Hence,
\begin{align*}
qN(c)-(q-1)^{m-1}&=qk_{m-2}(c^{1/m})-\sum_{x_1,\dots,x_{m-1}\in\F_q^*}\chi(x_1+\cdots+x_{m-1})\\
&=qk_{m-2}(c^{1/m})-(-1)^{m-1},
\end{align*}
and consequently
\[
	N(c)=k_{m-2}(c)+\frac1q((q-1)^{m-1}-(-1)^{m-1}),
\]
since $k_{m-2}(c^{1/m})=k_{m-2}(c)$. It now follows from~\eqref{e:N_p} that
\[
	N_m(a,b)=\frac{q^{m-1}-1}{q-1}+(-1)^{m-1}k_{m-2}(c),
\]
and the Deligne bound concludes the proof. 
\end{proof}

By Theorem~\ref{t:N_p_kl}, equation~\eqref{e:x_in_S_t}, and Lemma~\ref{l:polyN} we get an expression for $P_m(a,b)$ in terms of a hyper-Kloosterman sum:

\begin{corollary}\label{c:P_p} If $m=p^k>2$ and $ab\ne0$, then
\[
	mP_m(a,b)=\frac{q^{m-1}-1}{q-1}+(-1)^{m-1}k_{m-2}(b/a^m).
\]
\end{corollary}

\subsection{Irreducible cubics over $\F_{3^r}$}\label{cubics} 

Next we consider the number of irreducible cubics $P_3(a,b)$ when $q=3^r$. The main result of this section is the following:

\begin{corollary}Let $q=3^r$ and let $a,b\in\F_q$ with $ab\ne0$. Then, $P_3(a,b)=(q+1+t)/3$ where $t$ is an integer satisfying the following two conditions:
\begin{enumerate}
\item[(i)] $t\equiv-1\pmod 3$, 
\item[(ii)] $|t|<2\sqrt q$.
\end{enumerate}

Conversely, for a given integer $t$ satisfying conditions (i) and (ii) there are exactly 
$(q-1)H(t^2-4q)$ pairs $(a,b)\in\F_q^2$ with $ab\ne0$ and $P_3(a,b)=(q+1+t)/3$. Here $H(d)$ is the Kronecker class number of $d$.
\end{corollary}

\begin{proof}For a given $c\in\F_q^*$ there are exactly $q-1$ pairs $(a,b)\in\F_q^2$ such that $c=b/a^3$. Corollary~\ref{c:P_p}, Theorem~\ref{t:valuedist}, and Theorem~\ref{t:valuedist} below complete the proof.
\end{proof}

\begin{theorem}[{\cite{Katz-Livne}}]\label{t:valuedist}Let $q=3^r$. The range $S$ of $k_1(c)$, as $c$ runs over $\F_q^*$, is given by
\[
	S=\{t\in\Z : |t|< 2\sqrt q \text{ and } t\equiv-1\pmod 3\}.
\] 
Moreover, each value $t\in S$ is attained exactly $H(t^2-4q)$ times.
\end{theorem}
  
\begin{example}Let $q=3$. If $t$ is an integer satisfying conditions (i) and (ii) then $t=-1$ or $t=2$. There should be exactly $(3-1)H(1-12)=2$ pairs $(a,b)$ with $ab\ne0$ and $P_3(a,b)=1$, and exactly $(3-1)H(4-12)=2$ pairs $(a,b)$ with $ab\ne0$ and $P_3(a,b)=2$.  

Indeed, the two pairs $(a,b)$ for which there is exactly one irreducible polynomial $x^3+ax^2+cx+b\in\F_3[x]$ are $(a,b)=(1,1),(2,2)$, and the corresponding irreducible cubics are
\[
	x^3+x^2+2x+1,\ x^3+2x^2+2x+2. 
\] 
The two pairs $(a,b)$ for which there are exactly two irreducible cubics are $(a,b)=(1,2),(2,1)$ and the corresponding irreducible cubics are
\[
	x^3+x^2+2,\ x^3+x^2+x+2,\qquad x^3+2x^2+1,\ x^3+2x^2+x+1. 
\] 
Finally, for a pair $(0,b)$ there should be, by Example~\ref{ex:m_p_k}, exactly one irreducible cubic. Indeed, the corresponding polynomials are 
\[
	x^3+2x+1,\ x^3+2x+2. 
\]
Thus we have counted all the eight irreducible cubics in $\F_3[x]$.  
\end{example}

\subsection{Irreducible quartics over $\F_{2^r}$}\label{quartics}

We conclude Section~\ref{power} by considering the number of irreducible quartics $P_4(a,b)$ when $q=2^r$. We need the following result by Carlitz which links one- and two-dimensional Kloosterman sums:

\begin{theorem}[\cite{Carlitz69B}]\label{t:Carlitz} Let $c\in\F_q^*$. Then, 
\[
	k_2(c)=k_1(c)^2-q. 
\]
\end{theorem}

Now we are able to prove the main result of this section:

\begin{corollary}Let $q=2^r$ ($r>1$) and let $a,b\in\F_q$ with $ab\ne0$. Then, $P_4(a,b)=(q^2+2q+1-t^2)/4$ where $t$ is an integer satisfying the following two conditions:
\begin{enumerate}
\item[(i)] $t\equiv1\pmod 2$, 
\item[(ii)] $1\le t<2\sqrt q$.
\end{enumerate}

Conversely, for a given integer $\,t$ satisfying conditions (i) and (ii) there are exactly 
$(q-1)H(t^2-4q)$ pairs $(a,b)\in\F_q^2$ with $ab\ne0$ and $P_4(a,b)=(q^2+2q+1-t^2)/4$. 
\end{corollary}

\begin{proof} Let $c=b/a^4$. By Corollary~\ref{c:P_p} $4P_4(a,b)=q^2+q+1-k_2(c)$, and now, by Theorem~\ref{t:Carlitz}, we get
\[
	4P_4(a,b)=q^2+2q+1-k_1(c)^2.
\]
Theorem~\ref{t:valuedist2} below completes the proof.
\end{proof}

\begin{theorem}[\cite{La-Wo}]\label{t:valuedist2}Let $q=2^r$. The range $S$ of $k_1(c)$, as $c$ runs over $\F_q^*$, is given by
\[
S=\{t\in\Z : |t|< 2\sqrt q \text{ and } t\equiv-1\pmod 4\}.
\] 
Moreover, each value $t\in S$ is attained exactly $H(t^2-4q)$ times. 
\end{theorem}

\begin{example}Let $q=4$. Now $t=1$ and $t=3$ are the only integers satisfying (i) and (ii). There should be exactly $(4-1)H(1-16)=6$ pairs $(a,b)$ with $ab\ne0$ and $P_4(a,b)=(16+8+1-1)/4=6$,
and exactly $(4-1)H(9-16)=3$ pairs $(a,b)$ with $ab\ne0$ and $P_4(a,b)=(16+8+1-9)/4=4$.

Indeed, if $\F_4=\{0,1,\alpha,\beta\}$, then the six pairs $(a,b)$ for which there are exactly 
six irreducible polynomials $x^4+ax^3+\cdots+b\in\F_4[x]$ are 
\[
(a,b)=(1,\alpha),(1,\beta),(\alpha,1),(\alpha,\beta),(\beta,1),(\beta,\alpha),
\] 
and the three pairs $(a,b)$ for which there are exactly four irreducible quartics
\[
(a,b)=(1,1),(\alpha,\alpha),(\beta,\beta).
\]

Finally, for a pair $(0,b)$ there should be, by Example~\ref{ex:m_p_k}, exactly 
$(q^2+q+1-(q+1))/4=4$ irreducible quartics. This is indeed the case, and so we counted all the $6\cdot6+3\cdot4+4\cdot3=60$ irreducible quartics in $\F_4[x]$.

\end{example}

\section{Divisibility modulo three of Kloosterman sums, $q=2^r$}

Let $q=2^r$. We consider the divisibility modulo three of Kloosterman sums $k(c):=k_1(c)$. We use the following notations: 
\medskip

\begin{tabular}{cl}
$\Tr_{2^s}^q$ &\ the trace function from $\F_q$ onto $\F_{2^s}$\\
$A$ &\ the set of elements $a\in\F_q$ with $\Tr_{2}^q(a)=0$\\ 
$T_3(b)$ &\ the number of irreducibles $x^3+ax^2+cx+b\in\F_q[x]$\\
&\ with $b$ fixed and $a$ runs over the set $A$
\end{tabular}
\medskip

We need the following
\begin{theorem}[\cite{Moisio}]\label{t:Moisio} Let $\alpha\in\F_{q^m}^*$. Then,
\[
	\sum_{x\in\F_{q^m}^*}e(\alpha x^{q-1})=(-1)^{m-1}(q-1)k_{m-1}(\N_m(\alpha)). 
\]
\end{theorem}

\begin{lemma}\label{l:T_3}Let $b\in\F_q^*$. Then,
\[
	T_3(b)={\textstyle\frac13}\Big({\textstyle\frac12}(q^2+1+k(b)^2)-N(b)\Big),
\]
where $N(b)$ is equal to the number of solutions of $x^3=b$ in $A$.
\end{lemma}

\begin{proof} By Lemma~\ref{l:polyN}
\begin{equation}\label{e:T_3}
	3T_3(b)=\sum_{a\in A}N_3(a,b)-\sum_{a\in A}N_1(a,b).
\end{equation}
By~\eqref{e:x_in_S_t} the latter sum is equal to $N(b)$. Consider next the first sum. Apply Lemma~\ref{l:N_t} with $t=m=3$ to get
\[
	\sum_{a\in A}N_3(a,b)=\frac{q^3-1}{2(q-1)}+\frac1{q(q-1)}\sum_{a\in A}\sigma_3(a,b),
\]
where
\begin{align*}
\sum_{a\in A}\sigma_3(a,b)&=\sum_{c\in\F_q^*}\sum_{x\in\F_{q^3}^*}e(c\gamma_3^{i_0}x^{q-1})\sum_{a\in A}\chi(ca)\\
&=\frac12\sum_{c\in\F_q^*}\sum_{x\in\F_{q^3}^*}e(c\gamma_3^{i_0}x^{q-1})\sum_{a\in\F_q}\chi(c(a+a^2)).
\end{align*}
Since $\chi(ca)=\chi(c^2a^2)$ the orthogonality of characters now implies
\begin{align*}
\sum_{a\in A}\sigma_3(a,b)&=\frac12\sum_{c\in\F_q^*}
\sum_{x\in\F_{q^3}^*}e(c\gamma_3^{i_0}x^{q-1})\sum_{a\in\F_q}\chi((c+c^2)a^2)\\
&=\frac q2\sum_{x\in\F_{q^3}^*}e(\gamma_3^{i_0}x^{q-1}).	
\end{align*}
By Theorems~\ref{t:Moisio} and~\ref{t:Carlitz} we now get
\[
 \sum_{a\in A}\sigma_3(a,b)=\textstyle\frac12q(q-1)k_2(b)=\frac12q(q-1)(k_1(b)^2-q)
\]
(note that $\N_3(\gamma_3^{i_0})=g^{i_0}=b$), and therefore
\[
	\sum_{a\in A}N_3(a,b)=\textstyle\frac12(q^2+q+1+k_1(b)^2-q).
\]
Equation~\ref{e:T_3} now completes the proof.   
\end{proof}

\begin{theorem}\label{t:kl_mod3} Let $q=2^r$, and let $b\in\F_q^*$. Then, 3 divides $k(b)$ if and only if one of the following condition holds 
\begin{enumerate}
	\item $r$ is odd and $\Tr_2^q(\sqrt[3]{b})=0$
	\item $r$ is even, $b=a^3$ for some $a\in\F_q$, and $\Tr_4^q(a)\ne0$
\end{enumerate}
\end{theorem}  

\begin{proof} We have, by Lemma~\ref{l:T_3}
\[
	{\textstyle\frac12}(q^2+1+k(b)^2)-N(b)\equiv0\pmod 3,
\]
or, equivalently,
\[
k(b)^2\equiv -N(b)-2\pmod 3.
\]
Hence, $3\mid k(b)$ if and only if $N(b)\equiv 1\pmod 3$ if and only if $N(b)=1$. If $r$ is odd, then $x^3=b$ has unique solution $x=\sqrt[3]{b}$ in $\F_q$ and therefore $N(b)=1$ if and only if $\Tr_2^q(\sqrt[3]{b})=0$.

Assume $r$ is even and let $\zeta$ ($\in\F_4$) be a primitive third root of unity. Now $N(b)=1$ if and only if $b=a^3$ and $\Tr_2^q(a\zeta^i)=0$, for some $a\in\F_q$ and for unique $i\in\{0,1,2\}$. It follows by the transitivity of $\Tr_2^q$ that the latter condition is equivalent to $\Tr_4^q(a)\ne0$. 
\end{proof}

\begin{remark} In the case $r$ odd Theorem~\ref{t:kl_mod3} follows also from~\cite[Thm.~3]{chz} proved by using different methods.
\end{remark}

\section{A proof for the value distribution of a Kloosterman sum, $q=3^r$}

The aim of this section is to give a fairly elementary proof for Theorem~\ref{t:valuedist}. Let $q=3^r$, and let $c\in\F_q^*$.  Let $k(c):=k_1(c)$ and let $\mathcal X$ be the elliptic curve over $\F_q$ defined by 
\[
	\mathcal X : y^2+cy+xy=x^3.
\]  

\begin{lemma}\label{l:ellKl}
\[
	|\mathcal X(\F_q)|=q+1+k(c)\qquad\text{and}\qquad k(c)\equiv-1\pmod 3.
\]
\end{lemma}

\begin{proof}Choose $p=m=3$, and combine Theorem~\ref{t:P_3_a_ne0} and Corollary~\ref{c:P_p} to get
\[
	|\mathcal X(\F_q)|=3P(1,c)=q+1+k(c).
\]
\end{proof}

\begin{lemma}\label{l:ellisom} $\mathcal X$ is isomorphic over $\F_q$ to $\mathcal X':y^2=x^3+x^2-c$.
\end{lemma}

\begin{proof}Complete the square to get the equation of $\mathcal X$ in the form
\[
	(y+x+c)^2=x^3+(x+c)^2.
\] 
Then, substitute $y\mapsto y-x-c$, $x\mapsto x-c$ to get
\[
	y^2=x^3+x^2-c^3,
\]
and then, substitute $x\mapsto x^3$, $y\mapsto y^3$ to obtain
\[
	(y^2-x^3-x^2-c)^3=0.
\]
\end{proof}

Let $\mathcal E$ be an elliptic curve over $\F_q$. Starting from the long Weierstrass form 
\[
y^2+a_1xy+a_3y=x^3+a_2x^2+a_4x+a_6
\]
of the equation of $\mathcal E$, it is easy to see (see e.g.~\cite[p.10]{Washington}) that the equation of $\mathcal E$ can be given in the form
\[
y^2=x^3+ax^2+cx+d.
\]

If $a\ne0$ the substitution $x\mapsto x+e$ with $e=c/a$ yields the equation
\[
y^2=x^3+ax^2+e^3+ae^2+ce+d,
\]
and therefore we may assume that the equation of $\mathcal E$ is one of the following
\begin{enumerate}
\item[(i)] $y^2=x^3+ax^2+b$
\item[(ii)] $y^2=x^3+cx+b$
\end{enumerate}
for some $a,b,c\in\F_q$. Since $\mathcal E$ is smooth we must have $ab\ne0$ in case (i), and $c\ne0$ in case (ii). 

The $j$-invariant of $\mathcal E$ is given by
\[
j(\mathcal E)=
\begin{cases}
-a^3/b&\text{in case (i)},\\ 
0&\text{in case (ii)}.
\end{cases}
\]

\begin{lemma}\label{l:super}Let $|\mathcal E(\F_q)|=q+1+t$. The following three conditions are equivalent
\begin{enumerate}
\item $\mathcal E$ is supersingular,
\item $j(\mathcal E)=0$,
\item $3\mid t$.
\end{enumerate}
\end{lemma}

\begin{proof}See~\cite[p.75,p.121]{Washington}
\end{proof}

Assume next that $\mathcal E$ is ordinary (i.e. non-supersingular). We may now assume that $\mathcal E$ is defined by
\[
\mathcal E: y^2=x^3+ax^2+b.
\] 

\begin{lemma}\label{l:ordinary}If $a$ is a square in $\F_q^*$ then $\mathcal E$ is isomorphic over $\F_q$ to
\[
\mathcal X': y^2=x^3+x^2+b/a^3,
\]  
and $|\mathcal E(\F_q)|=q+1+t$ for some integer $\,t$ with $t\equiv-1\pmod 3$. 

If $a$ is not a square, then $|\mathcal E(\F_q)|=2(q+1)-|\mathcal X'(\F_q)|$, and 
$|\mathcal E(\F_q)|=q+1+t$ for some integer $t$ with $t\equiv1\pmod 3$.
\end{lemma}

\begin{proof}
If $a=c^2$ for some $c\in\F_q^*$, the substitutions $x\mapsto ax$, $y\mapsto c^3y$ yields the equation $y^2=x^3+x^2+b/a^3$. Assume next that $a$ is not a square. Let $\eta$ be the quadratic character of $\F_q$ with $\eta(0)=0$. The number of solutions $N$ of $y^2=x^3+ax^2+b$ in $\in\F_q^2$ is
\[
N=\sum_{x\in\F_q}(1+\eta(x^3+ax^2+b))=q+\sum_{x\in\F_q}\eta(x^3+ax^2+b).
\] 
Now substitute $x\mapsto ax$ to obtain
\[
N=q+\eta(a)\sum_{x\in\F_q}\eta(x^3+x^2+b/a^3)=q-\sum_{x\in\F_q}\eta(x^3+x^2+b/a^3),
\] 
and consequently
\[
|\mathcal E(\F_q)|=N+1=q+1-\big(|\mathcal X'(\F_q)|-(q+1)\big).
\] 
The remaining assertions follow now immediately by Lemmas~\ref{l:ellisom} and~\ref{l:ellKl}.
\end{proof}
 
\noindent{\em Proof of Theorem~\ref{t:valuedist}}. Let $t\equiv-1\pmod 3$ be an integer belonging to the interval $(-2\sqrt q,2\sqrt q)$. By Theorem~\ref{t:deuring} below there exist exactly $H(t^2-4q)$ pairwise non-isomorphic elliptic curves $\mathcal E$ with $|\mathcal E(\F_q)|=q+1+t$, and by Lemma~\ref{l:super} each of them is ordinary. Now, by Lemma~\ref{l:ordinary} each $\mathcal E$ is isomorphic over $\F_q$ to $\mathcal X':y^2=x^3+x+c$ for some $c\in\F_q^*$, and finally Lemmas~\ref{l:ellisom} and~\ref{l:ellKl} conclude the proof.\hfill\qed
  
\begin{theorem}[\cite{Deuring,Schoof}]\label{t:deuring}The number $M(t)$ of isomorphism classes of elliptic curves over $\F_q$ having $q+1+t$ points with $\gcd(q,t)=1$ is given by
\[
M(t)=
\begin{cases}
H(t^2-4q)&\text{if $t^2<4q$},\\
0&\text{otherwise}.
\end{cases}
\]
\end{theorem}

\begin{remark} Yet another proof of Theorem~\ref{t:valuedist}, which uses fairly advanced methods, is given in~\cite{Geer-Vlugt}.
\end{remark}

\end{document}